\documentclass[12pt,a4paper]{amsart}
\usepackage{amsmath,amsthm,amsfonts,amssymb,eucal, hyperref}

\newtheorem{theorem}{Theorem}[section]
\newtheorem{lemma}[theorem]{Lemma}
\newtheorem{proposition}[theorem]{Proposition}

\newtheorem{corollary}[theorem]{Corollary}

\newtheorem{remark}[theorem]{Remark}
\newtheorem{remarks}[theorem]{Remarks}

\allowdisplaybreaks
\renewcommand{\thefootnote}{\fnsymbol{footnote}}

\newcommand\finbox{~\hfill$\Box$}%

\def\la {{\lambda}}
\def\ep {{\epsilon}}

\newcommand {\pa}{\partial}

\numberwithin{equation}{section}
\begin{document}
{\centering
\bfseries
\title[Many Particle Hardy Inequalities]
{}
{\Large Many Particle Hardy Inequalities} 
\\[2\baselineskip]%
{\renewcommand\thefootnote{}%
\footnote{1991 Mathematics Subject Classification 35B05}}
\par
\mdseries
\scshape
\small
\author[Maria and Thomas Hoffmann-Ostenhof, A. Laptev,
J.Tidlom] {}
M. Hoffmann-Ostenhof$^1$\\
T. Hoffmann-Ostenhof$^{2,3}$ \\
A. Laptev$^4$\\
J. Tidblom$^5$
\par
\upshape
Fakult\"at f\"ur Mathematik , Universit\"at Wien$^1$\\
Institut f\"ur Theoretische Chemie, Universit\"at Wien$^2$\\
International Erwin Schr\"odinger Institute for Mathematical Physics$^3$\\
KTH, Stockholm$^4$\\
University of Stockholm$^5$\\
\today
}
\begin{abstract}
In this paper we prove three different types
of the so-called many-particle Hardy inequalities. 
One of them is a ``classical type" which is valid in any 
dimension $d\not=2$. The second type deals with two-dimensional
magnetic Dirichlet forms where every particle is supplied with a solenoid. 
Finally we show that Hardy inequalities for Fermions 
hold true in all dimensions.  

\end{abstract}
\maketitle


\section{Introduction}

Hardy inequalities play an important role in analysis.
The classical one states that for $u\in H_0^{1}(0,\infty)$
\begin{equation}\label{H1}
\int_0^\infty \Big|\frac{du}{dx}\Big|^2dx\ge \frac{1}{4}\int_0^\infty\frac{|u|^2}{|x|^2}dx.
\end{equation}
The standard Hardy inequality (away from a point) for functions  
$u\in H^{1}(\mathbb R^d)$ reads for $d\ge 3$
\begin{equation}\label{SHardy}
\int_{\mathbb R^d}|\nabla u|^2dx\ge \frac{(d-2)^2}{4}\int_{\mathbb R^d}\frac{|u|^2}{|x|^2}dx.
\end{equation}   
There are many other inequalities which also are called Hardy inequalities,
see for instance the survey paper by E.B. Davies \cite{D1} and the books 
of V.G. Maz'ya \cite{M} and Kufner and  Opic \cite{KO}.

In the present paper we shall investigate a kind of Hardy inequalities 
which might be called many-particle Hardy inequalities. They can 
be related to some Schr\"odinger operators 
and have  some interesting geometrical aspects.  

Pick $N$ a positive integer and consider $N$ particles. 
This means we consider $x\in \mathbb R^{dN}$, where $x=(x_1,x_2,\dots,x_N)$  
with $x_i=(x_{i,1},x_{i,2}\dots, x_{i,d})\in \mathbb R^d$. We define $r_{ij}$ by
\begin{equation*}\label{dist}
r_{ij}=|x_i- x_j|=\sqrt{\sum_{k=1}^d(x_{i,k}-x_{j,k})^2}.
\end{equation*}
We will write sometimes $\Delta_i=\sum_{k=1}^d\frac{\pa^2}{\pa x_{i,k}^2}$ so that 
$\Delta=\sum_{i=1}^N\Delta_i$. Similarily we write sometimes  
$\nabla_i$ for the gradient associated to the i-th particle. 

We have three groups of results. The first one deals with the  
``standard" Hardy inequality for many particles saying that 
\begin{equation}\label{sthardy}
 \sum_{j=1}^N\int_{\mathbb R^{dN}} |\nabla_{x_j} u|^2\,dx  \ge \,\mathcal  C(d,N)\,
\sum_{1\le i<j\le N} \int_{\mathbb R^{dN}} \frac{|u|^2}{r_{ij}^2}\,dx.
\end{equation}
In Sections \ref{4.1}-\ref{d=1} we prove that this inequality holds for 
$d\ge3$, $u\in H^1(\mathbb R^{dN})$ with a constant $\mathcal  C(d,N)$, such that 
$c_1N^{-1} \le \mathcal  C(d,N)\le c_2N^{-1}$, where $c_1,c_2>0$. 
The Hardy inequality \eqref{sthardy}
also holds for one-dimensional particles.
In this case the function $u$ is assumed to be equal to zero on diagonals
$x_i=x_j$. We find in this case that $\mathcal C(1,N)=1/2$ and 
that this constant is sharp.

In section \ref{magnetic} we consider the two-dimensional case 
and obtain a version of the Hardy inequality for magnetic multi-particle
Dirichlet forms with Aharonov--Bohm  type vector potentials attached 
to every particle. 
Let $x_j=(x_{j1},x_{j2})\in{\Bbb R}^2$, $j=1,2,\dots,N$, and let
\begin{equation}\label{ABfield}
{\bf F_j} = \alpha \Bigl( -\sum_{k\not= j} \frac{x_{j2}-x_{k2}}{r_{jk}^2},
\sum_{k\not= j} \frac{x_{j1}-x_{k1}}{r_{jk}^2}\Bigr),
\end{equation}
where $\alpha\in\Bbb R$. Then we shall prove that
\begin{equation}\label{magnhardy}
\int_{\mathbb R^{2N}}\sum_{j=1}^N|(i\nabla_{x_j} + {\bf F_j}) u|^2\,dx
\ge D_{N,\alpha}\, \int_{\mathbb R^{2N}}
|u|^2\Bigl(\sum_{k\not=j}\frac{1}{r_{kj}^2}\Bigr) \, dx.
\end{equation}
The explicit value for the constant $D_{N,\alpha}$ depends on
the "degree of rationality " of the magnetic flux $\alpha$.
 
Our third result concerns  the inequality \eqref{sthardy} for 
fermions, i.e. the anti-symmetric functions in $H^{1}(\mathbb R^{dN})$.
It turned out that in this case the Hardy 
inequality \eqref{sthardy} holds true in all dimensions and if $d\ge2$, then
\begin{equation}\label{fermions}
\mathcal  C(d,N)\ge \frac{d^2}{N},
\end{equation}
see Section \ref{fermions}.


\section{Main results}

\subsection{Hardy inequalities for $d$-dimenional particles with $d\ge 3$}

\begin{theorem}\label{d-dimH}
Assume that  $d\ge 3$, $N\ge 2$ and let $u\in H^{1}(\mathbb R^{dN})$.
Let us define
\begin{equation}\label{dNH}
\mathcal  C(d,N)=\inf_{u\in H^{1}(\mathbb R^{dN})}
\frac{\int_{\mathbb R^{dN}}|\nabla u|^2dx}{\sum_{1\le i<j\le N}
\int_{\mathbb R^{dN}}\frac{|u|^2}{r_{ij}^2}\;dx}.
\end{equation}
Then 
\begin{equation}\label{C+}
\mathcal  C(d,N)\ge(d-2)^2\max\Bigg\{\;\frac{1}{N},\; \frac{1}{1+\sqrt{1+
\frac{3(d-2)^2}{2(d-1)^2}(N-1)(N-2)}}\Bigg\}.
\end{equation}
\end{theorem}
%
\begin{remarks}\label{r2}

\noindent
\begin{enumerate}

\item
Hardy inequalities of this type cannot hold for general functions 
$u\in H^1(\Bbb R^{dN})$, $d=1,2$. 
\item
For large values of $N$ and $d\le6$ the maximum in \eqref{C+} is given 
by the second term.
\item
There is a very simple way of obtaining Hardy inequalities like above with 
a substantially weaker constant.  
Starting from \eqref{SHardy} and  noting that for any fixed $y\in \mathbb R^d$,
$d\ge3$, 
$$
-\Delta\ge \frac{(d-2)^2}{4}\frac{1}{|x-y|^2},
$$
we obtain 
$$
-\Delta_i-\Delta_j\ge  \frac{(d-2)^2}{2}\frac{1}{r_{ij}^2}
$$
in the quadratic form sense. Adding this up we would get 
\begin{equation*}\label{Drij}
-\Delta \ge \frac{(d-2)^2}{2N-2}\:\sum_{i<j}^N\frac{1}{r_{ij}^2}
\end{equation*}
in the sense of quadratic forms
and this is weaker than \eqref{C+} by a factor of more than 
 two for large $N$ and $d=3$.    
\item
The bounds for $\mathcal C(d,N)$ are not sharp. Actually for the lower bound
we use only the information from the derivation for the 3-particle case, 
i.e. $N=3$. There is certainly a lot of room for improvement, though it is not 
clear how to get explicit better bounds.
It is unclear what the optimal distribution of $\{x_j\}$ is as $N\to\infty$. 
Let $R(x,y,z)$ be the circumradius of the triangle with vertices $x$,$y$,$z$
and  suppose that the best asymptotic configuration of points 
could be described by a probability measure $\mu$ on ${\mathbb R}^d$. Let
\begin{equation*}
K= {\text{sup}}_\mu\,\,\frac{\int\int\int R^{-2}(x,y,z)\, \,d\mu(x)
d\mu(y)d\mu(z)}
{\int\int |x-y|^{-2}\,\,d\mu(x) d\mu(y)}\,.
\end{equation*}
Then applying \eqref{divmain}, see below, one can obtain a much better estimate 
of the constant $C(d,N)$ for large $N$
given by the inequality
$$
\lim_{N\to\infty} N\,\mathcal C(d,N)\ge \frac{(d-2)^2}{2+K}.
$$ 
Note that the integral 
$$
C^2(\mu) = \int\int\int R^{-2}(x,y,z)\, \,d\mu(x)
d\mu(y)d\mu(z)
$$
is known as {\it Menger-Melnikov curvature of the measure $\mu$},
see \cite{MMV}, \cite{Tol}.  Finding the value of $K$ is an interesting 
open problem.

\end{enumerate}
\end{remarks}



The next theorem shows that the estimate $\mathcal  C(d,N) = O(N^{-1})$, as $N\to\infty$,
cannot be improved. 
\begin{theorem}\label{C<<}
Let $\varphi\in H^{1}(\mathbb R^d)$, $d\ge3$, and define 
\begin{equation*}\label{conv}
M(\varphi)=\int_{\mathbb R^d\times \mathbb R^d}\frac{|\varphi(x)|^2|\varphi(y)|^2}
{|x-y|^2}dxdy
\end{equation*}
and 
\begin{equation*}
\mathcal D(d)=\inf_{\varphi\in C_0^\infty(\mathbb R^d)}\:
\frac{\int_{\mathbb R^d}|\nabla \varphi|^2dx \int_{\mathbb R^d}|\varphi|^2dx}
{M(\varphi)}.
\end{equation*}
Then  
\begin{equation}\label{Dd}
\mathcal C(d,N)\le \frac{2\mathcal D(d)}{N-1}.
\end{equation} 
\end{theorem}

For numerical upper bounds see \eqref{C<} and also Remark \ref{><}.

\begin{corollary}\label{Unb}
For any $N$ and $d\ge 3$ there is a constant ${\mathcal C}'={\mathcal C}'(d,N)$ 
such that the operator in $\Bbb R^{dN}$
$$
-\Delta-{\mathcal C}'(d,N)\sum_{1\le i<j\le N}\frac{1}{r_{ij}^2}
$$
is not bounded from below and such that 
$$
{\mathcal C}'(d,N) \le \frac{c(d)}{N}
$$ 
with some $c(d)>0$.
\end{corollary}

Corollary \ref{Unb} could be obtained by explicit calculation, see \eqref{Tphi} 
and \eqref{rijphi}. Indeed, 
it follows from Theorem \ref{C<<} that there is a function, 
$\varphi\in H^{1}(\mathbb R^{dN})$ and 
a constant $c$ such that 
\begin{equation}\label{phi-b}
\int_{\mathbb R^{dN}}|\nabla \varphi|^2-\frac{c}{N}\,\sum_{1\le i<j\le N}
\int_{\mathbb R^{dN}}
\frac{1}{r_{i,j}^2}|\varphi|^2dx<0.
\end{equation}
Since  both $-\Delta$ and the potential term show the same scaling, then if we
replace $\varphi(x)$ by $\varphi(\la x)$ and normalize we can make the expression in \eqref{phi-b} as 
negative as we want.


\subsection{Hardy inequality for 1D particles}

\begin{theorem}\label{1dHar}
Let 
\begin{equation*}
\mathcal N_N=\{x=(x_1,x_2,\dots,x_N)\in{\mathbb R}^N\:\big|\: x_i=x_j\,\,
{\text{for some}}\,\, i\neq j\}.
\end{equation*}
Suppose that $u\in H_0^{1}({\mathbb R}^N\setminus\mathcal N_N)$ then
\begin{equation}\label{1dN}
\int_{\mathbb R^N}|\nabla u|^2dx\ge \frac{1}{2}\int_{\mathbb R^N}|u|^2\sum_{i<j}^N\frac{1}
{r_{ij}^2}dx.
\end{equation}
The constant $1/2$ is sharp.
\end{theorem}

\begin{remark}\label{1drem}
One can get easily an  inequality like \eqref{1dN} with an N-dependent constant
instead of $1/2$ by using \eqref{H1}. First note that \eqref{H1} can be rewritten
such that for any $y\in{\mathbb R}$ and $u\in H_0^{1}({\mathbb R}\setminus\{y\})$ 
\begin{equation}\label{1dHy}
\int_{-\infty}^\infty \Big|\frac{du}{dx}\Big|^2dx\ge \frac{1}{4}\int_{-\infty}^\infty
\frac{|u|^2}{|x-y|^2}dx.
\end{equation}
Now consider $N=2$ and note that \eqref{1dHy} implies for 
$u\in H_0^{1}({\mathbb R}^2\setminus{\mathcal N}_2)$ 
\begin{equation*}
\int_{\mathbb R_2}\Big|\frac{\partial u}{\partial x_i}\Big|^2dx
\ge \frac{1}{4}\int_{\mathbb R^2}\frac{|u|^2}{|x_1-x_2|^2}dx
\end{equation*}
for $i=1,2$ so that adding up we get $\|\nabla u\|^2\ge \frac{1}{2}\|u/r_{12}\|^2$.
If we would continue like this we would get instead of the $1/2$ in \eqref{1dN},
$\frac{1}{2N-2}$ as in \eqref{Drij} of Remark \ref{r2}, 
a much weaker bound tending to zero for 
$N\rightarrow \infty$. 
\end{remark}


\subsection{Magnetic Hardy inequalities in 2D}

Let the vector field ${\bf F} = ({\bf F}_1,{\bf F}_2, \dots, {\bf F}_N)$ be defined by 
\eqref{ABfield} and let 
\begin{equation}\label{D_N}
D_{N,\alpha}=\min_{l=1,\dots,N-1}\Bigl(\frac{\min_{k\in{\Bbb Z}}|k-l\alpha|}{l}\Bigr)^2.
\end{equation}
\begin{theorem}\label{Magn}
{\footnote[1]{The energy integral in the right 
hand side of \eqref{magnhardy} appears when studying the Fractional 
Quantum Hall Effect. It has been considered in \cite{Laugh}, \cite{Jain}, where the fractional filling factor has been explained by 
attaching to each electron an infinitely thin magnetic 
solenoid carrying an Aharonov-Bohm flux (each electron 
bound to a flux tube has been called a ``composite particle", \cite{Jain}). 
}}
The following magnetic Hardy inequality for two-dimensional particles 
holds true
\begin{equation*}
\int_{\mathbb R^{2N}}\sum_{j=1}^N|(i\nabla_{x_j} + {\bf F_j}) u|^2\,dx
\ge D_{N,\alpha}\, \int_{\mathbb R^{2N}}
|u|^2\Bigl(\sum_{k\not=j}\frac{1}{r_{kj}^2}\Bigr) \, dx.
\end{equation*}
\end{theorem}

This inequality could be considered as a version of a 2D Hardy 
inequality by Laptev-Weidl \cite{LW} for Aharonov-Bohm magnetic 
Dirichlet forms and its generalisation obtained by 
A. Balinsky \cite{Bal}.


\subsection{Hardy inequalities for fermions}

Let us consider anti-symmetric functions $u(x)$, $x=(x_1, x_2, \dots, x_N)$, 
$x_j\in {\mathbb R}^d$,  such that
$$
u(x_1,\dots, x_i,\dots, x_j,\dots, x_N)=-u(x_1,\dots, x_j,\dots, x_i, \dots, x_N)
$$
 for all pairs  $(i,j),\; i\neq j$.
 
\begin{theorem}\label{Ferm} 
For any  $d=1,2, \dots$, and anti-symmetric function $u\in H^1({\mathbb R^{dN}})$
we have
 \begin{equation}\label{FermHardy}
 \sum_j \int_{{\mathbb R}^{dN}} |\nabla_{x_j} u|^2\, dx \ge
\frac{d^2}{N}\,
 \sum_{i<j} \int_{{\mathbb R}^{dN}}\frac{ | u|^2}{r_{ij}^2}\, dx.
 \end{equation}
 \end{theorem}

 \begin{remark}
The latter inequality could be improved for large $N$. 
By using arguments from \cite{LThirring} and  \cite{LYau} we expect that for large $N$ the
$N$ dependence of the constant in \eqref{FermHardy}  could be improved to $N^{-1/3}$.

It has recently been 
shown in \cite{FHLS} that there is a constant $C_d$ such that 
$\mathcal  C(d,N)\ge C_d\, N^{-1/3}$.
\end{remark}

\section{Some auxiliary results}

In this section we consider several simple results of analytical and 
geometrical character and start with a simple but crucial inequality. 
%
\begin{lemma}\label{div}
Let $u\in H^{1}(\mathbb R^m), \;m\ge 1$ and let 
\begin{equation*}\label{fatF}
\mathbf{\mathcal F}=\big(\mathcal F_1(x),\mathcal F_2(x), \dots, \mathcal F_m(x)\big)
\end{equation*}
be a  vectorfield in $\mathcal F\,:\mathbb R^m\mapsto \mathbb R^m$ whose components  and their first derivatives are uniformly bounded in
$\mathbb R^m$.   Then 
\begin{equation}\label{I}
\int_{\mathbb R^m}|\nabla u|^2dx\ge \frac{1}{4}\,
\frac{\Big(\int_{\mathbb R^m}|u|^2{\rm{div}}\, \mathbf{\mathcal F}dx\Big)^2}
{\int_{\mathbb R^m}|u|^2|\mathbf{\mathcal F}|^2dx}.
\end{equation}
\end{lemma}


\noindent
\textbf{Proof.}\\
We use the Cauchy-Schwarz inequality and partial integration. Indeed,
\begin{multline*}\label{I2}
\Big|\int_{\mathbb R^m}|u|^2{\rm{div}}{\mathcal F}dx\Big|=2|{\Re} \int_{\mathbb R^m}
\langle {\mathcal F},\, \nabla u\rangle \overline u dx|\\
\le 2\Big(\int_{\mathbb R^m}|u|^2|\mathbf{\mathcal F}|^2dx\Big)^{1/2}
\Big(\int_{\mathbb R^m}|\nabla u|^2dx\Big)^{1/2}.  
\end{multline*}
Squaring this inequality completes the proof.
\finbox

\noindent
The standard Hardy inequality (away from a point),
\eqref{SHardy} for $m\ge 3$
can be obtained by choosing 
\begin{equation*}
\mathcal F_\ep=\frac{x}{|x|^2+\ep^2}.
\end{equation*}
We pick  $u\in H^1(\mathbb R^m)$ and insert $\mathcal F_\ep$ into \eqref{I}
and obtain 
\begin{equation*} 
{\rm div}{\mathcal F}_\ep=\frac{m}{|x|^2+\ep^2}-\frac{2|x|^2}{(|x|^2+\ep^2)^2}
\ge \frac{m-2}{|x|^2+\ep^2}
\end{equation*}
so that 
\begin{equation*}
\int_{\mathbb R^m}|\nabla u|^2dx\ge \frac{(m-2)^2}{4}\,
\int_{\mathbb R^m}\frac{|u|^2}{|x|^2+\ep^2}dx.
\end{equation*}
  $C_0^\infty(\mathbb R^d)$ is dense in $H^1(\mathbb R^d)$ and therefore $\ep\rightarrow 0$ 
gives \eqref{SHardy}.


The next lemma is related to the so-called Melnikov-Menger 
curvature and could be found,  for example, in \cite{Tol}. 
\begin{lemma}\label{circum} 
Define for three points $x_i,x_j,x_k\in \mathbb R^d$,  
\begin{equation*}\label{bijk}
b_{ijk}=\frac{\langle x_i-x_j,\; x_i-x_k\rangle}{r_{ij}^2r_{ik}^2}
+\frac{\langle x_j-x_i,
\;x_j-x_k\rangle}{r_{ij}^2r_{jk}^2}
+\frac{\langle x_k-x_i,\; x_k-x_j\rangle}{r_{ik}^2r_{jk}^2}.
\end{equation*}
 Let $R_{ijk}$ be the circumradius
of the triangle with corners $x_i,x_j,x_k$. Then 
\begin{equation*}\label{cradius}
b_{ijk}=\frac{1}{2R_{ijk}^2}, \quad {\rm{if}} \quad d\ge2 \qquad {\rm{and}}
\qquad b_{ijk} =0, \quad  {\rm{if}} \quad d=1.
\end{equation*}
\end{lemma}

\noindent
\textbf{Proof.} Let $a=x_i-x_j$ and $b=x_i-x_k$. 
Then 
\begin{multline}\label{r}
b_{ijk}=\frac{(a\cdot b)}{|a|^2|b|^2}-
\frac{a\cdot(b-a)}{|a|^2|b-a|^2}-\frac{b\cdot(a-b)}{|b|^2|b-a|^2}\\
=\frac{2\big(|a|^2|b|^2-(a\cdot b)^2\big)}{|a|^2|b|^2|b-a|^2}=
\frac{2\sin^2\phi}{r_{jk}^2}.
\end{multline}
Here $\phi$ is the angle between $a$ and $b$. The relation between 
the circumradius and the angle follws from the sine-theorem. Clearly if 
$x_i,x_j,x_k\in \Bbb R$ and not all of them equal, $R_{ijk} = \infty$.
\finbox\\

The next statements are concerned with two inequalities for triangles, see also 
\cite{Mitr}. 
\begin{lemma}\label{tri}
Let $R$ be the circumradius of a triangle with sides with side lenghts
$a,b,c$ then 
\begin{equation}\label{trtr}
\frac{1}{R^2}\le \frac{9}{a^2+b^2+c^2}\le \frac{1}{a^2}+
\frac{1}{b^2}+\frac{1}{c^2}.
\end{equation}
Both inequalities are equalities for the  equilateral triangle.
\end{lemma}

\noindent
\textbf{Proof.}
This is an easy consequence of the sine-theorem 
and a Lagrange multiplyer argument. 
Indeed notice that $R=\frac{a}{2\sin \alpha}=\frac{b}{2\sin \beta}=\frac{c}{2\sin \gamma}$ where the 
angles $\alpha, \beta, \gamma$ correspond to the angle at the corner opposite to the sides with 
side lengths  $a$ respectively $b, c$. We show the first inequality. This reads 
\begin{equation}\label{9R2}
4R^2(\sin^2\alpha+\sin^2\beta+\sin^2\gamma)\le 9R^2. 
\end{equation}
It hence suffices to show that for 
$\alpha+\beta+\gamma=\pi$, $\sin^2\alpha+\sin^2\beta
+\sin^2\gamma\le 9/4$. 
So we look at 
\begin{equation*}\label{max!}
\sin^2\alpha+\sin^2\beta+\sin^2\gamma +\la(\alpha+\beta+\gamma-\pi)=\max\:!   
\end{equation*}
Differentiation leads to 
\begin{equation*}
\sin 2\alpha+\la =0, \;\sin 2\beta+\la =0,\; \sin 2\gamma +\la=0, \; \alpha+\beta+\gamma=\pi
\end{equation*}
and this implies that 
\begin{equation*}\label{sinsin}
\sin 2\alpha=\sin 2\beta =\sin 2\gamma,\:\: \alpha+\beta+\gamma =\pi.
\end{equation*}
There are three solutions, namely $\alpha=\beta=\gamma$, $\alpha=\beta=\pi/2,\: \gamma=0$ and finally
$\alpha =\pi, \;\beta=\gamma=0$. If we insert the values into \eqref{9R2} we get the desired result.

For the other inequality we have to show that  
\begin{equation*}
P=\Big(\frac{1}{a^2}+\frac{1}{b^2}+\frac{1}{c^2}\Big)(a^2+b^2+c^2)\ge 9.
\end{equation*}
and this can be seen by multiplication which yields
$$
P=3+a^2/b^2+b^2/a^2+a^2/c^2+c^2/a^2+b^2/c^2+c^2/b^2\ge 9
$$
since $a^2/b^2+b^2/a^2\ge 2$ and similarily for the other fractions above.

\finbox

The following two statements can be checked by straight forward 
computations.
\begin{lemma}\label{MM}
 Let $x_j\in {\mathbb R}^d$, $j=1,2,3$. Then
\begin{multline*} 
r_{12}^2 + r_{13}^2 + r_{23}^2\\
 =2\,\Bigl(\langle x_1 - x_2, x_1 - x_3\rangle +  \langle x_2 - x_1, x_2 - x_3\rangle 
+\langle x_3 - x_1, x_3 - x_2\rangle\Bigr).   
 \end{multline*}
 \end{lemma}
\begin{lemma}\label{NN}
Let  $x_j\in {\mathbb R}^d$, $j=1,2,\dots, N$. Then
\begin{equation*} 
N\, \sum_{j=1}^N \Delta_{x_j} \\
 = \sum_{1\le j<k\le N} (\nabla_{x_j} -\nabla_{x_k})^2 
 + \Bigl(\sum_{j=1}^N  \nabla_{x_j} \Bigr)^2.
\end{equation*}
 \end{lemma}
%
 
Finally we need a statement which could be considered 
as two versions of Hardy's inequalities for three particles.

\begin{lemma}\label{rho2}
Let $x_1,x_2,x_3\in \mathbb R^d$, $d\ge 2$, and let
\begin{equation*}\label{ro2}
\rho^2=r_{12}^2+r_{13}^2+r_{23}^2.
\end{equation*}
Then  
\begin{equation}\label{Hrh}
\int_{\mathbb R^{3d}}|\nabla u|^2dx\ge 3(d-1)^2\int_{\mathbb R^{3d}}\frac{|u|^2}{\rho^2}dx.
\end{equation}
Furthermore if $R(x)$ is the circumradius of the triangle with vertices $x_1, x_2, x_3$,
then
\begin{equation}\label{RH}
\int_{\mathbb R^{3d}}|\nabla u|^2dx\ge \frac{(d-1)^2}{3}\int_{\mathbb R^{3d}}\frac{|u|^2}{R^2}dx.
\end{equation}
\end{lemma} 

\noindent
\textbf{Proof.}
This follows from a simple direct calculation.  
Let $\mathcal F=\mathcal G$ in \eqref{I}, where
\begin{equation}\label{Grho}
 \mathcal G=\frac{1}{\rho^2}\Big( 2x_1-x_2-x_3,\:2x_2-x_1-x_3,\:2x_3-x_1-x_2\Big).
\end{equation}
Then by applying Lemma \ref{div}
we easily work out by using the identity  given in Lemma \ref{MM}, that
\begin{equation*}\label{rhodeG}
{{\rm div}\, \mathcal G}=\frac{6(d-1)}{\rho^2}
\end{equation*}
and 
\begin{equation*}\label{rhoG2}
|\mathcal G|^2=\frac{3}{\rho^2}.
\end{equation*}
We insert these equalities
into \eqref{I} and obtain \eqref{Hrh}. To be more precise we first consider 
$\mathcal G_\ep$ where  the denominator in \eqref{Grho} is replaced by $\rho^2+\ep^2$.
Then as in the proof of the standard Hardy-inequality the result follows as $\ep$ tends
to zero.  
Finally in order to prove \eqref{RH} we use the inequality from Lemma \ref{tri}, 
which tells us that 
\begin{equation*}\label{9R}
\rho^2\le 9R^2.
\end{equation*}
Hence \eqref{RH} follows immediately 
from \eqref{Hrh}.\finbox

\begin{remarks}

\noindent
\begin{enumerate}
\item
For one-dimensional particles the circumradius is equal to infinity and 
therefore \eqref{RH} becomes trivial. 
\item
However, we do not believe that the constant in \eqref{RH} is sharp. 
Perhaps one can find a suitable $\mathcal F$ so that 
one can directly obtain a Hardy-type inequality for $R^{-2}$. 
\end{enumerate} 
\end{remarks}

\section{Proofs of main results.}

\subsection{Proof of Theorem \ref{d-dimH}}\label{4.1}

\subsection*{A}
Let us first give a simple proof of the inequality \eqref{dNH} which states that
$\mathcal  C(d,N)\ge(d-2)^2/N$.

For a function $u\in H^1(\Bbb R^{dN})$ we consider a vector field 
$$
{\mathcal F_1}(x_j,x_k) = (x_j-x_k)r_{jk}^{-2}, \qquad 1\le j<k\le N.
$$
Then by using arguments from the proof of Lemma \ref{div} with the vector field 
$\mathcal F_1$ we find
\begin{multline}\label{jk}
\int_{\Bbb R^{dN}} |(\nabla_{x_j} - \nabla_{x_k}) u|^2\, dx_j dx_k\\
\ge \frac{1}{4} 
\frac{\Big(\int_{\mathbb R^{dN}}|u|^2\,\bigl(({\rm div}_{x_j} - {\rm div}_{x_k})\, 
\mathbf{\mathcal F_1}\bigr)\, dx\Big)^2}
{\int_{\mathbb R^{dN}}|u|^2|\mathbf{\mathcal F_1}|^2dx}\\
= (d-2)^2 \int_{\Bbb R^{dN}} \frac{|u|^2}{r_{jk}^2}\, dx_j dx_k.
\end{multline} 
Moreover, if we introduce the vector field 
$$
\mathcal F_2(x) = \frac{\sum_{j=1}^N  x_j}{\Big| \sum_{j=1}^N  x_j \Big|^2},  
$$
then using  Lemma \ref{div} with $\mathcal F_2$ we obtain 
\begin{equation}\label{dN}
\int_{\Bbb R^{dN}} \Big|\sum_{j=1}^N \nabla_{x_j} u\Big|^2 
\ge \frac{(d - 2)^2N^2}{4} \int_{\Bbb R^{dN}} 
\frac{|u|^2}{\Big|\sum_{j=1}^N x_j \Big|^2} \, dx.
\end{equation}
Adding the inequalities \eqref{jk} and \eqref{dN} up and using Lemma \ref{NN} 
we arrive at
\begin{multline}
\int |\nabla u|^2 dx \ge 
\frac{(d-2)^2}{N}\, \int_{\Bbb R^{dN}} \sum_{j<k} \frac{|u|^2}{r_{jk}^2}\, dx\\
+ \frac{(d - 2)^2N^2}{4N} \int_{\Bbb R^{dN}} 
\frac{|u|^2}{\Big|\sum_{j=1}^N x_j \Big|^2}\,dx.
\end{multline}
The latter inequality implies the inequality $\mathcal  C(d,N)\ge(d-2)^2/N$ and also 
gives a positive remainder term which is of order $O(N)$.


\subsection*{B}

Let us now define
\begin{equation}\label{F_3}
 \mathcal F_3=(F_1,\dots, F_N),
\end{equation}
where the $F_j$ are given by 
\begin{equation}\label{F_3j}
F_j=\sum_{k\neq j}^N\frac{x_j-x_k}{r_{jk}^2}.
\end{equation}
In order to prove Theorem \ref{d-dimH}
we apply Lemma \ref{div} for the
vector field ${\mathcal F_3}$ which is 
conveniently written as a vector with $N$ elements which themselves are 
vectors with $d$ entries. The divergence of $\mathcal F_3$ can be similarily defined as 
\begin{equation*}\label{dNdiv}
{\rm{div}}\, \mathcal F_3=\sum_{i=1}^N\nabla_i\cdot F_i
\end{equation*}
where $\nabla_i\cdot F_i=$ div $F_i$ and where the divergence is now with 
respect to a $d$-dimensional vector field.

\begin{proposition}\label{XR}
Assume that $d\ge 3$ and $N\ge 2$. Let for an arbitrary $u\in H^{1}(\mathbb R^{dN})$ 
\begin{gather*}\label{defX}
T(d,N)=\int_{\mathbb R^{dN}}|\nabla u|^2dx,\\ 
X(d,N)=\sum_{1\le 1<j\le N}\int_{\mathbb R^{dN}}
\frac{1}{r_{ij}^2}|u|^2dx
\end{gather*}
 and 
\begin{gather*}
Z(d,N)=\sum_{1\le i<j<k\le N}\int_{\mathbb R^{dN}}\frac{1}{R_{ijk}^2}|u|^2dx,
\end{gather*}
where   $R_{ijk}$ is as in Lemma \ref{circum}.  Then
\begin{equation}\label{divmain}
T(d,N)\ge (d-2)^2\:\frac{X(d,N)^2}{2X(d,N)+Z(d,N)}.
\end{equation}
\end{proposition}

\noindent
\textbf{Proof.}\\
The proof is an easy calculation. We just note that  
\begin{equation}\label{divF_3}
{\rm{div}\;\mathcal F_3}=2(d-2)\sum_{1\le i<j\le N}\frac{1}{r_{ij}^2}
\end{equation}
and that 
\begin{equation}\label{F2_3}
|\mathcal F_3|^2=2\sum_{1\le i<j\le N}\frac{1}{r_{ij}^2}+\sum_{1\le i<j<k\le N}
\frac{1}{R_{ijk}^2},
\end{equation}
where we used Lemma \ref{circum}.
We just have to insert these expressions into \eqref{I} to obtain 
\eqref{divmain} proving the proposition. \finbox

Consider now inequality \eqref{divmain}. There are two possibilities to 
obtain from this quadratic inequality a linear inequality

{\bf a.} First we can try to find an estimate such that 
\begin{equation*}\label{KdN}
Z(d,N)\le k(d,N)X(d,N)
\end{equation*}
and this leads to 
\begin{equation}\label{casea}
\mathcal C(d,N) \ge \frac{(d-2)^2}{2+k(d,N)}.
\end{equation}
{\bf b.} The other possibility is to find an estimate of the form
\begin{equation*}\label{LdN}
Z(d,N)\le \ell(d,N)T(d,N).
\end{equation*}
Indeed, with this estimate we get 
\begin{equation*}\label{sqL}
T(d,N)\ge(d-2)^2\frac{X(d,N)^2}{2X(d,N)+\ell(d,N)T(d,N)}
\end{equation*}
and this leads to the quadratic inequality
\begin{equation}\label{Q<}
X(d,N)^2-\frac{2X(d,N)T(d,N)}{(d-2)^2}-\frac{\ell(d,N)T(d,N)^2}{(d-2)^2}\le 0.
\end{equation}
Therefrom we get by solving the corresponding quadratic equation
\begin{equation}\label{Xsq}
\mathcal \mathcal  C(d,N)\ge \frac{(d-2)^2}{1+\sqrt{1+\ell(d,N)^2(d-2)^2}}.
\end{equation}
\textbf{case a.}\\
We show that 
\begin{equation}\label{C-N}
\mathcal  C(d,N)\ge \frac{(d-2)^2}{N}.
\end{equation}
This is an easy consequence of the inequality \eqref{trtr} in Lemma \ref{tri}. Indeed
we just have to show  \eqref{casea} that $k(d,N)\le N-2$ and this can be seen
by counting. Clearly $Z(d,N)$ consists of $\binom{N}{3}$ and $X(d,N)$ 
of $\binom{N}{2}$ terms. Finally
we group each three particle coordinates together and apply Lemma \ref{tri}.
This gives \eqref{C-N} and an alternative proof of the result obtained in 
subsection {\bf A}.

\noindent
\textbf{case b.}\\
This case is more involved. We begin with considering three particles. 

Note that  for three $d$-dimensional particles with $d\ge 3$, \eqref{RH}
implies
\begin{equation*}\label{ell<}
\ell(d,3)\le \frac{3}{(d-1)^2}
\end{equation*}
so that we obtain for $N=3$ in \eqref{Xsq}
\begin{equation*}\label{Xsq3}
T(d,3)\ge \frac{(d-2)^2}{1+\sqrt{1+\frac{3(d-2)^2}{d-1)^2}}}X(d,3).
\end{equation*}
We continue with the $N$-particle case and get by counting from \eqref{divmain}
that
\begin{equation*}\label{RT}
T(d,N)\ge \frac{2(d-1)^2}{3(N-1)(N-2)}\;Z(d,N).
\end{equation*}
From the quadratic inequality \eqref{Q<} we now infer that 
\begin{equation*}\label{NHa2}
T(d,N)\ge (d-2)^2\frac{1}{1+\sqrt{1+\frac{3(d-2)^2}{2(d-1)^2}\:(N-1)(N-2)}}
X(d,N).
\end{equation*}
This inequality together with \eqref{C-N} proves \eqref{C+} and therefore the second 
part of Theorem \ref{d-dimH}.
\finbox

\subsection{Proof of Theorem \ref{C<<}.}\label{UpperBound}
Let $\varphi(x_i)\in H^{1}(\mathbb R^d)$ and consider for fixed $N$ 
\begin{equation*}\label{PhiN}
u(x)=u_N(x)=\prod_{i=1}^N\varphi(x_i),\:\: x_i\in \mathbb R^d.
\end{equation*}
We observe that 
\begin{equation}\label{Tphi}
\int_{\mathbb R^{dN}}|\nabla u_N(x)|^2dx=
N\int_{\mathbb R^d}|\nabla_1\varphi(x)|^2dx\Bigg(\int_{\mathbb R^d}|\varphi(x)|^2dx\Bigg)^{N-1}.
\end{equation}
Next we calculate 
\begin{multline}\label{rijphi}
\sum_{1\le i<j\le N}\int_{\mathbb R^{dN}}\frac{1}{r_{ij}^2}|u_N|^2dx=\\
\frac{N(N-1)}{2}
\int_{\mathbb R^d}\int_{\mathbb R^d}|\varphi(x)|^2\frac{1}{|x-y|^2}|
\varphi(y)|^2dxdy\:\Bigg(\int_{\mathbb R^d}|\varphi(x)|^2dx\Bigg)^{N-2}.
\end{multline}
Substituting the expressions from \eqref{Tphi} and \eqref{rijphi} into
\eqref{dNH} we complete the proof.
\finbox

\bigskip
Here we provide a numerical value for the right hand side in \eqref{Dd}
and therefore an estimate from above for the constant ${\mathcal C}(d,N)$.

Let us choose 
$$
\varphi(x) = e^{-|x|^2/2}. 
$$
Then 
$$
\int_{{\Bbb R}^d} |\varphi(x)|^2\, dx 
= \frac{1}{2}\, |{\Bbb S}^{d-1}|\, \Gamma(d/2), \qquad
\int_{{\Bbb R}^d} |\nabla \varphi(x)|^2\, dx
=\frac{d}{4}\, |{\Bbb S}^{d-1}|\,  \Gamma(d/2).
$$
Straight forward computations give us
\begin{equation*}
\int_{{\Bbb R}^d\times{\Bbb R}^d} 
\frac{|\varphi(x)|^2\,|\varphi(y)|^2}{|x-y|^2}\, dx dy = 
\frac{2\pi^{d/2}}{\Gamma(d/2)}.
\end{equation*}
Substituting all the expressions into \eqref{Dd} we obtain
\begin{equation}\label{C<}
{\mathcal C}(d,N) \le  \frac{2d}{2(N-1)}\, \pi^{d/2}\, \Gamma(d/2).
\end{equation}
In particular,  
$$
{\mathcal C}(3,N) \le \frac{1}{N-1} \, \frac{3\pi^2}{4},
$$
i.e. $0.43<{\mathcal C}(3,N)<3.69$.
For the three particle system using the estimate from below 
provided by Theorem \ref{d-dimH} we have
$$ 
\frac{1}{1+\sqrt7/2} \le {\mathcal C}(3,3) \le \frac{3\pi^2}{8}.
$$

\begin{remark}\label{><}
It follows from \eqref{C+} that the gap between the lower 
and upper bounds obtained in Theorem \ref{d-dimH} and  in formula 
\eqref{C<} is growing with respect to $d$.
\end{remark} 

\subsection{Proof of Theorem \ref{1dHar}.}\label{d=1}

%

The inequality \eqref{1dN}
 follows immediately from Lemma \ref{div} with $\mathcal F$ defined by 
 \eqref{F_3}, \eqref{F_3j} and the relations
\eqref{divF_3} and \eqref{F2_3}. It only remains to observe that by Lemma \ref{circum}
the second sum in \eqref{F2_3} is equal to zero for $d=1$. 

Let us now prove that the constant $1/2$ appearing in \eqref{1dN}
is sharp. It is enough to show that  for any $\varepsilon >0$ there is a 
function $v=v_\varepsilon$ such that 
\begin{equation}\label{1/2function}
\int_{\Bbb R^N} |\nabla v(x)|^2\, dx \le \Big(\frac12 + \varepsilon\Big) \, 
\int_{\Bbb R^N} |v(x)|^2 \sum_{i<j}^N \frac{1}{r_{ij}^2}\, dx.
\end{equation}
Let $\alpha = 1/4 + \delta$ 
\begin{equation*}
v(x) = \Pi_{i\not=j} (x_i-x_j)^{2\alpha} e^{-|x|}.
\end{equation*}
Then 
$$
\partial_{x_i} v = 2\alpha v\,\sum_{j:j\not=i}\frac{1}{x_i-x_j} - v \, \frac{x_i}{|x|} .
$$
Therefore 
\begin{multline}\label{nabla v}
|\nabla v|^2 = \sum_{i=1}^N |\partial_{x_i} v|^2\\ 
= \sum_{i=1}^N \Big(4\alpha^2\, v^2
\Big(\sum_{j:j\not=i}\frac{1}{(x_i-x_j)^2}
+\sum_{j,k:j,k\not=i,\, j\not=k} \frac{1}{x_i-x_j}\, \frac{1}{x_i-x_k}\Big)\\
-4\alpha v^2 \,\sum_{j:j\not=i}\frac{1}{x_i-x_j}\,\frac{x_i}{|x|} 
+ v^2\,\frac{x_i^2}{|x|^2}\Big).
\end{multline}
Note that by Lemma \ref{circum}
$$
\sum_{i=1}^N\sum_{j,k:j,k\not=i,\, j\not=k} \frac{1}{x_i-x_j}\, \frac{1}{x_i-x_k}
= \sum_{i=1}^N\sum_{j,k:j,k\not=i,\, j\not=k} \frac{x_i-x_j}{r_{ij}^2}\, \frac{x_i-x_k}{r_{ik}^2} =0.
$$
Moreover the identity 
$$
\sum_{j:j\not=i}\frac{1}{x_i-x_j}\,\frac{x_i}{|x|} = \frac{1}{|x|}\, \sum_{j:j\not=i}
\Big(1 - \frac{x_j}{x_j-x_i}\Big)
$$
implies 
$$ 
\sum_{i=1}^N\sum_{j:j\not=i}\frac{1}{x_i-x_j}\,\frac{x_i}{|x|} = \frac{1}{2|x|}\, N(N-1)\ge0.
$$
Therefore we obtain from \eqref{nabla v} 
\begin{multline*}
\int_{\Bbb R^N} |\nabla v|^2\, dx = \int_{\Bbb R^N} \Big(\sum_{i<j}
\frac{8\alpha^2\, v^2}{r_{ij}^2} - \frac{2\alpha v^2}{|x|}\, N(N-1)+v^2\Big)\, dx\\
\le 8\alpha^2 \int_{\Bbb R^N} \sum_{i<j}
\frac{v^2}{r_{ij}^2}\, dx\,  \Big(1 + \beta(\delta)\Big),
\end{multline*}
where 
$$
\beta(\delta) = \frac{\int_{\Bbb R^N} v^2\,dx}{8(1/4+\delta)^2\,\int_{\Bbb R^N} 
v^2\,\sum_{i<j}\frac{1}{r_{ij}^2}\,dx} \to 0 \quad {\rm{as}}\quad \delta\to 0.
$$
We conclude the proof by choosing $\delta$ 
small enough so that it satisfies the inequality
$1/2+\varepsilon \ge 8(1/4+\delta)^2(1+\beta(\delta))$.
\finbox

\subsection{Proof of Theorem \ref{Magn}.}\label{magnetic}

We begin with recalling two results obtained in the papers
of \cite{LW} and \cite{Bal} concerning the Hardy inequalities  
for Aharonov-Bohm magnetic Dirichlet forms.  

\noindent
{\bf a. One particle inequality.} 
Let $x=(x_1,x_2)\in\mathbb R^2$, $\alpha\in{\mathbb R}$ and let ${\bf{F}}$ be the  
Aharonov-Bohm vector potential
$$
{\bf{F}} = (F_1, F_2) = \alpha\,\Bigl(-\frac{x_2}{|x|^2}, \frac{x_1}{|x|^2}\Bigr).
$$
\begin{lemma} 
$$
\int_{\mathbb R^2} |(i\nabla + {\bf F})u|^2\, dx 
\ge \min_{k\in{\mathbb Z}} \,(k-\alpha)^2\, \int_{\mathbb R^2} \frac{|u|^2}{|x|^2}\, dx.
$$
\end{lemma}

\noindent
{\bf Proof.}
Indeed, using polar coordinates $(r,\theta)$ we have 
$u(x)=\frac{1}{\sqrt{2\pi}}\sum_k u_k(r)e^{ik\theta}$.
Therefore
\begin{multline*}
\int_{\mathbb R^2} |(i\nabla + {\bf F})u|^2\, dx 
= \int_0^\infty\int_0^{2\pi}\Bigl(|u'_r|^2 
+ \Big|\frac{iu'_\theta + \alpha u}{r}\Big|^2\Bigr)\,r\,d\theta\,dr\\
\ge \frac{1}{2\pi}\, \int_0^\infty\int_0^{2\pi}\Big|\sum_k 
\frac{\alpha - k}{r}\,u_k e^{ik\theta}\Big|^2 \,r\,d\theta\,dr
= \int_0^\infty  \sum_k\Big|\frac{\alpha - k}{r}\,u_k\Big|^2 \,r\,d\theta\,dr\\
\ge \min_{k\in{\mathbb Z}} \,(k-\alpha)^2\, \int_{\mathbb R^2} \frac{|u|^2}{|x|^2}\, dx.
\end{multline*}
\finbox

\bigskip
\bigskip
\bigskip
\bigskip
\bigskip
\bigskip

\noindent
{\bf b. Magnetic potentials with multiple singularities.}
 
Assume that  $\{z_1, z_2, \dots, z_n\}$ are $n$ fixed different 
points in ${\Bbb C}$, $z_j = x_j+iy_j$ and 
$\alpha_j\in{\Bbb R} $. Let ${\bf F}$ the following vector potential 
$$
{\bf F} = \sum_{j=1}^n \alpha_j \, 
\Bigl(\frac{-y+y_j}{|z-z_j|^2}, \frac{x-x_j}{|z-z_j|^2}\Bigr), \quad z=x+iy.
$$
This corresponds to Aharonov-Bohm magnetic vector fields placed 
in $n$ points $z_j$ with magnetic fluxes $\alpha_j$.
Let now $\Phi: {\Bbb C} \to {\Bbb C}$ be an analytic function with zero set 
 $\{z_1, z_2, \dots, z_n\}$ and such that $\Phi(\infty)=\infty$. 

\noindent
Let $\{\xi_1, \xi_2, \dots, \xi_m\}$ be the zero set of $\Phi_z'$
and let  $\{0, |\Phi(\xi_1)|, \dots, |\Phi(\xi_m)|\}$ be such that 
$0\ge |\Phi(\xi_1)|\ge \dots \ge |\Phi(\xi_m)|$.
Denote by  $\mathcal A$ the 
pre-image of  these points under the map $|\Phi|:\, {\Bbb C} \to {\Bbb R}_+$. 
For an arbitrary point $z\not\in\mathcal A$ we define a curve $\gamma_z$
obtained by $|\Phi|^{-1}(|\Phi(z)|)$. 
Let $\Omega_z\subset {\Bbb C}$
be a bounded  domain defined  by $\gamma_z$. 
We now consider a piecewise constant function
\begin{equation}\label{C_Phi}
C_\Phi(z) = \Bigl(\frac{\text{min}_{k\in{\Bbb Z}} \, |k- \sum_{j:\, z_j\in\Omega_z} \alpha_j|}{\sum_{j:\, z_j\in\Omega_z} 1}\Bigr)^2.
\end{equation}
\begin{lemma} {\rm (A. Balinsky)}
The following Hardy inequality holds true
$$
\int_{{\Bbb R}^2} |(i\nabla + {\bf F}) u|^2 \, dxdy 
\ge \int_{{\Bbb R}^2} C_\Phi(z) \Big|\frac{\Phi_z'}{\Phi}\Big|^2 \, |u|^2\,dx dy.
$$
\end{lemma}
For the proof see \cite{Bal}.

\noindent
{\bf c. Multi-particle case.}
Let now $z=(z_1,\dots,z_N)$,  $z_j=x_{j1} + ix_{j2}$ 
and let $\Phi_j(z) = \Pi_{k\not=j}(z_j-z_k)$, 
$j,k=1,\dots,N$. Then according Balinsky's lemma there are piecewise constants functions $C_{\Phi_j}(x)$ defined by \eqref{C_Phi}, such that
$$
\int_{\mathbb R^{2N}}|i\nabla_{x_j} + {\bf F_j}) u|^2\,dx
\ge \int_{\mathbb R^{2N}} C_{\Phi_j}(x)\Big|\frac{(\Phi_j)_{z_j}'(z)}{\Phi_j(z)}\Big|^2\,|u|^2\,dx.
$$
A simple computation shows
$$
\Big|\frac{(\Phi_j)_{z_j}'(z)}{\Phi_j(z)}\Big|^2 
= \Big|\sum_{k\not=j} \frac{1}{z_j-z_k}\,\Big|^2
=\sum_{k,l\not=j} \frac{(x_j-x_k)\cdot(x_j-x_l)}{r_{jk}^2 r_{jl}^2}\,.
$$
Note that $C_{\Phi_j}(x)\ge D_{N,\alpha}$, where 
$D_{N,\alpha}$ is defined by \eqref{D_N}.
Therefore we obtain
\begin{multline*}
\int_{\mathbb R^{2N}}\sum_{j=1}^N|(i\nabla_{x_j} + {\bf F_j})u|^2\,dx
\ge D_{N,\alpha}\, \int_{\mathbb R^{2N}}
\sum_{j=1}^N \Big|\sum_{k\not=j}^N \frac{1}{z_j-z_k}\Big|^2\,|u|^2\,dx \\
= D_{N,\alpha}\, \int \Bigl(\sum_{k\not=j}^N \frac{1}{r_{jk}^2} 
+ \sum_{l\not=k, l,k\not=j}^N\frac{1}{R_{jkl}^2}\Bigr)|u|^2\, dx.
\end{multline*}
We complete the proof by noticing that  
$$
\min_{x\in{\Bbb R}^{2N}} \sum_{l\not=k, l,k\not=j}^N R_{jkl}^{-2}=0.
$$
\finbox

\subsection{Proof of Theorem \ref{Ferm}.}\label{fermions}

Let us begin with a simple observation concerning odd functions 
in $\mathbb R^d$ which has been pointed out already in the classical paper 
of M.S. Birman \cite{Bir}. 
\begin{proposition}\label{odd}
Let $u(x)=-u(-x)\in H^1(\mathbb R^d)$, $d\ge 2$. Then 
\begin{equation*}
\int_{\mathbb R^d}|\nabla u|^2dx \ge \frac{d^2}{4}
\int_{\mathbb R^d}\frac{|u|^2}{|x|^2}dx.
\end{equation*}
\end{proposition}

\noindent
\textbf{Proof.}\\
Let us introduce spherical coordinates $x = (r,\theta)$. Then 
$$
\int_{\mathbb R^d}|\nabla u|^2dx
=\int_0^\infty\int_{\mathbb S^{d-1}} \Bigl(|u_r'|^2 
+ \frac{|\nabla_\theta u|^2}{r^2}\Bigr)\, r^{d-1}\, d\theta dr.
$$
By using the 1-dimensional Hardy inequality with weight (see for example \cite{KO}) 
we obtain
\begin{multline*}
\int_0^\infty\int_{\mathbb S^{d-1}} |u_r'|^2 r^{d-1}\, d\theta dr
\ge \frac{(d-2)^2}{4}\,   \int_0^\infty\int_{\mathbb S^{d-1}} |u|^2 r^{d-3}\, d\theta dr\\
= \frac{(d-2)^2}{4}\,
\int_{\mathbb R^d}\frac{|u|^2}{|x|^2}dx.
\end{multline*}
It only remains to note that since $u$ is an odd function
it is orthogonal to constants on $\mathbb S^{d-1}$ and therefore 
\begin{multline*}
\int_0^\infty\int_{\mathbb S^{d-1}} |\nabla_\theta u|^2\, r^{d-3}\, d\theta dr
\ge (d-1) \int_0^\infty\int_{\mathbb S^{d-1}} |u|^2\, r^{d-3}\, d\theta dr\\
=(d-1)\,\int_{\mathbb R^d}\frac{|u|^2}{|x|^2}dx,
\end{multline*}
where $d-1$ is the second eigenvalue of the Laplace-Beltrami 
operator on $\mathbb S^{d-1}$.
\finbox

\bigskip
\noindent
We now consider an anti-symmetric function of two variables
$x,y\in{\mathbb R}^d$.
 \begin{lemma}\label{anti-symm}
 For any anti-symmetric function $u(x,y)=-u(y,x) \in H^1(\mathbb R^{2d})$
 we have
\begin{equation*}
\int_{\mathbb R^{2d}}|(\nabla_x-\nabla_y) u(x,y)|^2dxdy\ge 
d^2 \int_{\mathbb R^4}\frac{|u(x,y)|^2}{|x-y|^2}dxdy.
\end{equation*}
\end{lemma}

 \noindent
 {\bf Proof.}
 We make an orthogonal  coordinate transformation 
\begin{equation*}
s=\frac{1}{\sqrt 2}(x+y),\:\: t=\frac{1}{\sqrt 2}(x-y).
\end{equation*}
Thus $|x|^2+|y|^2=|s|^2+|t|^2$ and 
$$
\nabla_s=\frac{1}{\sqrt 2}(\nabla_x+\nabla_y),\:\:\nabla_t=
\frac{1}{\sqrt 2}(\nabla_x-\nabla_y), \:\:
\Delta=\Delta_x+\Delta_y=\Delta_s+\Delta_t.
$$ 
If we define the function $\tilde u(s,t)$ as
$$
\tilde u(s,t)= u(x,y)  = u\Big(\frac{s+t}{\sqrt2},\frac{s-t}{\sqrt2}\Big),
$$
then it is odd with respect to $t$, $\tilde u(s,-t) = -\tilde u(s,t)$.
By using Proposition \ref{odd} we obtain
\begin{equation*}
\int_{\mathbb R^{2d}}|\nabla_t\tilde u(s,t)|^2\, dsdt \ge 
 \frac{d^2}{4}\int_{\mathbb R^{2d}}
\frac{|\tilde u|^2}{|t|^2}dsdt.
\end{equation*}
 Transforming back to $u$ and noting that $|t|^{-2}=2|x-y|^{-2}$ 
 we complete the proof.
\finbox

\noindent
Let us note (cf. Lemma \ref{NN})  that for $\xi\in \Bbb R^{dN}$ 
 \begin{equation*}
 \sum_{j=1}^N |\xi_j|^2 = 
\frac{1}{N}\, \sum_{j<k}  |\xi_j - \xi_k |^2 
+ \frac{1}{N}\, \Big|\sum_{j=1}^N \xi_j \Big|^2.
 \end{equation*}
\noindent
If $\hat u$ is the Fourier transform of the function $u$,
then by using Lemma \ref{anti-symm} we find
\begin{multline*}
 \sum_j \int_{{\mathbb R}^{dN}} |\nabla_{x_j} u|^2\, dx
 =  \sum_j \int_{{\mathbb R}^{dN}} |\xi_j \hat u|^2\, d\xi \\
 =  \frac{1}{N}\,
\sum_{j<k} \int |(\xi_j - \xi_k)\hat u |^2\, d\xi +
 \frac{1}{N}\,
\int |\sum_j  \xi_j\,\hat u |^2\, d\xi\\
 =\frac{1}{N}\,
\sum_{j<k} \int |(\nabla_{x_j} - \nabla{x_k})\, u |^2\, dx 
+ \frac{1}{N}\,
\int \Big|\sum_{j=1}^N \nabla_{x_j}  u \Big|^2\, dx \\
 \ge 
\frac{d^2}{N}\, \sum_{i<j} \int_{{\mathbb R}^{dN}}\frac{| u|^2}{r_{ij}^2}\, dx.
\end{multline*}
 In the latter inequality we can neglect the second integral and this completes 
 the proof of Theorem \ref{Ferm}.
\finbox 

\medskip
\noindent
{\it Acknowledgements.}
The authors are grateful to partial support by the ESF European Programme  "SPECT". 
A. Laptev and J. Tidblom would like to thank the International Erwin Schr\"odinger Institute for its hospitality. T. Hoffmann-Ostenhof has been partially supported by G\"oran Gustafsson Foundation. The first three authors would like to express their gratitude for the hospitality 
of the International Newton Institute.

\bibliographystyle{amsplain}

\scshape 
M. Hoffmann-Ostenhof: Fakult\"at f\"ur Mathematik, Universit\"at Wien,
Nordbergstrasse 15, A-1090 Wien, Austria.

email: Maria.Hoffmann-Ostenhof@univie.ac.at

\scshape
T. Hoffmann-Ostenhof: Institut f\"ur Theoretische Chemie, Universit\"at
Wien, W\"ahringer Strasse 17, A-1090 Wien, Austria and International Erwin
Schr\"odinger Institute for Mathematical Physics, Boltzmanngasse 9, A-1090
Wien, Austria.

email: thoffman@esi.ac.at

\scshape
A. Laptev: Department of Mathematics, KTH, 100 44 Stockholm, Sweden.

email: laptev@math.kth.se

\scshape 
J. Tidblom: Department of Mathematics, Univerity of Stockholm, 
106 91 Stockholm, Sweden.

email: jespert@math.su.se

\end{document}